\documentclass[12pt]{article}
\usepackage{amssymb,latexsym,amsmath}

\newcommand\Res{\operatorname{Res}} 
\newcommand\res{\operatorname{res}}

\makeatletter
\newfont{\footsc}{cmcsc10 at 8truept}
\newfont{\footbf}{cmbx10 at 8truept}
\newfont{\footrm}{cmr10 at 10truept}
\renewcommand{\ps@plain}{%
\renewcommand{\@oddfoot}{\footsc the electronic journal of combinatorics {\footbf 8} (2001), \#N7\hfil\footrm\thepage}}
\makeatother
\title{The polynomial part of a restricted partition function related to the Frobenius problem}

\author{ Matthias Beck\\
\small Department of Mathematical Sciences \\ [-0.8ex]
\small State University of New York \\[-0.8ex]
\small Binghamton, NY 13902--6000, USA \\ [-0.8ex]
\small \texttt{matthias@math.binghamton.edu}
\and Ira M.~Gessel
\thanks{Research partially supported by 
NSF grant DMS-9972648. } \\ 
\small Department of Mathematics \\  [-0.8ex]
\small Brandeis University \\  [-0.8ex]
\small Waltham, MA 02454--9110, USA \\  [-0.8ex]
\small \texttt{ gessel@brandeis.edu}
\and
\ Takao Komatsu \\ 
\small Faculty of Education \\  [-0.8ex]
\small Mie University \\  [-0.8ex]
\small Mie, 514--8507, Japan \\  [-0.8ex]
\small \texttt{ komatsu@edu.mie-u.ac.jp} }

\date{\small Submitted: May 29, 2001; Accepted: September 4, 2001 \\
\small MR Subject Classifications: Primary 05A15; Secondary 11P81, 05A17}

\begin{document}
\maketitle

\begin{abstract}
\noindent Given a set of positive integers $ A = \{ a_{1} , \dots , a_{n} \} $, 
we study the number $ p_{A} (t) $ of nonnegative integer solutions 
$ \left( m_{1} , \dots , m_{n} \right) $ to $ \sum_{j=1}^{n} m_{j} a_{j} = t $. 
We derive an explicit formula for the polynomial part of $p_A$.
\end{abstract}

\bigskip
Let $ A = \{ a_{1} , \dots , a_{n} \} $ be a set of positive integers 
with $\gcd(a_{1}, \dots, a_{n}) = 1$. The classical \emph{Frobenius problem} 
asks for the largest integer $t$ (the \emph{Frobenius number}) such that 
  \[ m_1 a_1 + \dots + m_n a_n = t \]  
has no solution in nonnegative integers $ m_{ 1 } , \dots , m_{ n } $. 
For $n=2$, the Frobenius number is $(a_1-1)(a_2-1)-1$, as is well known,
but the problem is extremely difficult for
$n>2$.  (For surveys of the Frobenius problem, see \cite{alfonsin,selmer}.) 
One approach \cite{bdr,israilov,komatsu,sertoz} is to study  the restricted partition function $ p_{A} (t)
$, the number of nonnegative integer solutions 
$ \left( m_{1} , \dots , m_{n} \right) $ to $ \sum_{j=1}^{n} m_{j} a_{j} = t $, 
 where $t$ is a nonnegative integer. The Frobenius number 
is the largest integral zero of $ p_{A} (t) $. Note that, in contrast  to the 
Frobenius problem, in the definition of $p_A$ we do not require 
$a_{1}, \dots, a_{n}$ to be relatively prime. 
In the following, $a_{1}, \dots, a_{n}$ are \emph{arbitrary} positive integers. 
 
It is clear that $ p_A (t) $ is the coefficient of
$z^t$ in the generating function 
  \[ G(z)=\frac 1 { ( 1 - z^{ a_1 } ) \cdots ( 1 - z^{ a_n } ) } . \] 
If we expand $G(z)$ by partial fractions, we see that
$p_A(t)$ can be written in the form 
\[ \sum_\lambda P_{A,\lambda}(t) \lambda^t,\]
where the sum is over all complex numbers $\lambda$ such that $\lambda^{a_i}=1$ for some $i$, 
and $P_{A,\lambda}(t)$ is a polynomial in $t$.
The aim of this paper is to give an explicit formula for $ P_{A,1}(t)$, which we denote by $P_A(t)$ and call the
\emph{polynomial part} of $p_A(t)$. It is easy to see that $P_A(t)$ is a polynomial of degree $n-1$.
(More generally, the degree of $P_{A,\lambda}(t)$ is one less than the number of values of $i$ for which
$\lambda^{a_i}=1$.) It is well known \cite[Problem 27]{polya} that 
\[ p_{A} (t) = \frac{ t^{n-1} }{ (n-1)! \, a_{1} \cdots a_{n} } + O \left( t^{ n-2 } \right) . \] 
Our theorem is a refinement of this statement. 
We note that Israilov derived a more complicated formula for $P_A(t)$ in \cite{israilov}. 

Let us define $Q_A(t)$ by $p_A(t)=P_A(t)+Q_A(t)$.
From the partial fraction expansion above, it is clear that $Q_A$ (and 
hence also $p_A$) is a \emph{quasi-polynomial}, that is, an expression of the form
  \[ c_{d}(t) t^{d} + \dots + c_{1}(t) t + c_{0}(t), \]
where $ c_{0}, \dots , c_{d} $ are periodic functions in $t$. 
(See, for example, Stanley \cite[Section 4.4]{stanley}, for more information about quasi-polynomials.)
In  the special case in which the $a_i$ are pairwise relatively prime,
each $P_{A,\lambda}(t)$ for $\lambda\ne1$ is a constant, and thus $Q_A(t)$ is a periodic
function with average value 0, and this property determines $Q_A(t)$, and thus $P_A(t)$.
Discussions of $Q_A(t)$ can be found, for example, in  \cite{bdr,israilov,komatsu}. 

We define the \emph{Bernoulli numbers} $B_j$ by 
  \begin{equation}\label{bernoulli} \frac z {e^z - 1} = \sum_{j \ge0}  B_j\frac{ z^j}{ j! } \end{equation}  
(so $ B_0 = 1 , B_1 = - \frac 1 2 , B_2 = \frac 1 6 , B_4 = - \frac 1 {30}$,
and $B_n=0$ if $n$ is odd and greater than 1.)

{\medskip\noindent\bf Theorem.} 

\begin{align}
 P_A (t) &= \frac 1 { a_1 \cdots a_n } \sum_{ m=0 }^{ n-1 } \frac{ (-1)^m }{ (n-1-m)! } \sum_{  k_1 + \dots + k_n = m  } a_1^{ k_1 }  \cdots a_n^{ k_n } \frac{ B_{ k_1 } \cdots B_{ k_n } }{ k_1! \cdots k_n! } t^{n-1-m}\label{e:t1}\\
         &=\frac 1 { a_1 \cdots a_n } \sum_{ m=0 }^{ n-1 } \frac{ (-1)^m }{ (n-1-m)! }\notag\\
         &\kern 20pt\times\kern -5pt\sum_{  k_1 + 2 k_2 \dots + m k_m = m } \frac{(-1)^{k_2+\cdots+k_m} }{ k_1! \cdots k_m! } \left(\frac{ B_1 s_1}{1\cdot 1!} \right)^{ k_1 } \cdots \left( \frac{ B_m s_m }{ m \cdot m! } \right)^{ k_m } t^{n-1-m},\label{e:t2}
\end{align} 
{\it where} $s_i=a_1^i+\cdots+a_n^i$.

As an intermediate step, we first prove the following more elegant but 
less explicit formula for $P_A(t)$. 

{\medskip\noindent\bf Proposition.} {\it $P_A(t)$ is the constant term in $z$ in} 
  \[ - \frac{ z e^{ -tz } }{ ( 1 - e^{ a_1 z } ) \cdots ( 1 - e^{ a_n z } ) } \ . \] 

{\smallskip\noindent\it Proof.} As noted earlier, $ p_A (t) $ is the coefficient of $z^t$ in
the generating function 
  \[ G(z)=\frac 1 { ( 1 - z^{ a_1 } ) \cdots ( 1 - z^{ a_n } ) } . \] 
Hence if we let $f(z)=G(z)/z^{t+1}$
then $ p_A (t) = \Res ( f(z) , z=0 ) $. 
As in \cite{bdr}, we use the residue theorem to derive a formula for $ p_A (t) $. 
Since clearly $ \displaystyle \lim_{R \to \infty} \int_{|z|=R} f(z) \, dz = 0 $, 
  \[ p_A (t) = - \Res ( f(z) , z=1 ) - \sum \Res ( f(z) , z=\lambda ) . \] 
Here the sum is over all nontrivial $a_1, \dots , a_n$th roots of unity $\lambda$.  
It is not hard to see that
$\Res(f(z), z=\lambda)$
may be expressed in the form $u_\lambda(t)\lambda^{-t}$ for some polynomial $u_\lambda(t)$,
and thus it follows from our earlier discussion that $-\Res(f(z), z=\lambda)=P_{A,\lambda^{-1}}(t)\lambda^{-t}$.
In particular,
  \[ P_A (t) = - \Res ( f(z) , z=1 ) . \] 
To compute this residue, note that 
  \[ \Res ( f(z) , z=1 ) = \Res ( e^z f(e^z) , z=0 ), \] 
so that 
  \begin{equation}\label{res} P_A (t) = - \Res \left( \frac{ e^{ -tz } }{ ( 1 - e^{ a_1 z } ) \cdots ( 1 - e^{ a_n z } ) } , z=0 \right). \end{equation}  

\hfill {} $\Box$ 

{\smallskip\noindent\it Proof of the theorem.} 
The coefficient of $t^{n-1-m}$ in $ P_A (t) $ is by (\ref{res}) the coefficient of $z^{ -n+m }$ in 
\[ \frac{ (-1)^{n-m} }{ (n-1-m)! }\cdot\frac 1 {(1-e^{a_1 z})\cdots(1-e^{a_n z})},\]
which is the coefficient of $z^m$ in 
\begin{equation}\frac{(-1)^m}{(n-1-m)!\, a_1\cdots a_n} B(a_1z)\cdots B(a_nz),
\label{e:bprod}
\end{equation}
where $B(z)=z/(e^z-1)$, and this implies \eqref{e:t1}.

To prove \eqref{e:t2}, we first note that 
\[\log\left(\frac{z}{e^z-1}\right)=\sum_{j\ge1}(-1)^{j-1} \frac{B_j}j\frac{z^j}{j!},\]
as can easily be proved by differentiating both sides.
Then 
\begin{align}
 B(a_1z)\cdots B(a_nz) &= \exp \sum_{j\ge 1} (-1)^{j-1}\frac{  B_j s_j }{j} \frac{ z^j}{j!}\notag\\
                       &=\prod_{j\ge1} \exp\left((-1)^{j-1}\frac{  B_j s_j }{ j} \frac{ z^j}{j!} \right). \label{e:bpowersum}
\end{align}
Since $B_{2i+1}=0$ for $i>0$, $(-1)^{j-1}B_j=-B_j$ for $j>1$,  
and \eqref{e:t2} follows from \eqref{e:bprod} and \eqref{e:bpowersum}.

\hfill {} $\Box$ 

{\noindent\it Remark.} It is possible to avoid the use of complex analysis and give a purely formal power series
proof of the theorem. We indicate here how this can be done. We work with formal Laurent series, which are power
series with finitely many negative powers of the variable. If
$F(z)=\sum_i u_i z^i$ is a formal Laurent series ($u_i$ is nonzero for
only finitely many negative values of $i$) then the {\it residue\/} of
$F(z)$ is $\res F(z)=u_{-1}$. An elementary fact about formal Laurent series is the change
of variables formula for  residues: If $g(z)$ is a formal power series
with $g(0)=0$ and $g'(0)\ne0$ then
\[\res F(z)=\res F(g(z)) g'(z).\]
(See, for example, Goulden and Jackson \cite[p.~15]{gj}.)

\def\N{G}
\def\F{U}
By partial fractions, we have 
\[\N(z)=\frac{1}{ (1-z^a_1)\cdots (1-z^{a_m})}=\frac{c_1}{ 1-z}+\cdots+\frac{c_m}{(1-z)^m}+R(z),\] 
where $R(z)$ is a rational function of $z$ with denominator not divisible by $1-z$.
It follows from our earlier discussion that
\[\sum_{t=0}^\infty P_A(t) z^t=\frac{c_1}{ 1-z}+\cdots+\frac{c_m}{(1-z)^m}\]
and thus
\[P_A(t)= \sum_{l=1}^\infty c_l {t+l-1\choose l-1},\] 
where we take $c_l$ to be 0 for $l>m$.
Now let $\F(z)=\N(1-z)$. Then
\begin{align*}\F(z)&= \frac{1}{(1-(1-z)^{a_1})\cdots(1-(1-z)^{a_m})} \\
                   &=\frac{c_1}{z} +\cdots+\frac{c_m}{z^m} +R(1-z),\end{align*} 
where $R(1-z)$ has a formal power series expansion (with no negative powers
of $z$), and thus $c_l=\res z^{l-1}\F(z)$. Note that this holds for all $l\ge1$,
since for $l>m$, $c_l=\res z^{l-1}\F(z)=0$.

Then 
\[P_A(t)=\sum_{l=1}^\infty c_l {t+l-1\choose l-1}=\res\frac{\F(z)}{z}\sum_{l=1}^m z^l \binom{t+l-1}{ l-1}=\res \frac{\F(z)}{(1-z)^{t+1}}.\]

\def\P{a_1\cdots a_m}
\def\B{B(a_1z)\cdots B(a_mz)}
\def\P{a_1\cdots a_m}

We now apply the change of variables formula with $g(z)=1-e^z$ and we obtain
\begin{align*}
P_A(t)&=-\res {\F(1-e^z)\over e^{tz}}\\
      &=-\res \frac{e^{-tz}} {(1-e^{a_1z})\cdots (1-e^{a_mz})},
\end{align*}
which is \eqref{res}, and the proof continues as before.


\small 
\addcontentsline{toc}{subsubsection}{References}
\bibliography{thesis}
\bibliographystyle{alpha}

\end{document}